\documentclass[preprint,12pt,letterpaper,reqno]{amsart}
\RequirePackage{amsthm,amsmath,amsfonts,amssymb}
\usepackage[only,llbracket,rrbracket]{stmaryrd}
\usepackage{dsfont}
\usepackage{xcolor}
\usepackage{mathrsfs}
\usepackage{enumerate}
\usepackage{float}
\RequirePackage[colorlinks,citecolor=green!40!black!80,urlcolor=orange]{hyperref}
\RequirePackage{graphicx}
\def\R{\mathds{R}}
\def\N{\mathds{N}}

\allowdisplaybreaks

\DeclareMathOperator{\supp}{Supp}
\DeclareMathOperator{\diam}{diam}

\newcommand{\me}{\mathrm{m}}

\newcommand{\G}{\mathscr{G}}
\numberwithin{equation}{section}
\theoremstyle{plain}
\newtheorem{thm}{Theorem}[section]

\newtheorem{coro}{Corollary}[section]
\newtheorem{lema}{Lemma}[section]

\theoremstyle{remark}

\newtheorem{dfn}{Definition}[section]

	\author[Castillo]{Ren\'e Erlin Castillo$^{1}$}
	\address[$^{1}$]{Departamento de Matem\'aticas\newline\indent 
		Universidad Nacional de Colombia,\newline\indent 
		Bogot\'a, Colombia.}%
	\email{recastillo@unal.edu.co}
\author[Chaparro]{H\'ector Camilo Chaparro$^{2}$}
\address[$^{2}$]{Programa de Matem\'aticas\newline\indent 
	Universidad de Cartagena,\newline\indent 
	Cartagena de Indias, Colombia.}%
\email{hchaparrog@unicartagena.edu.co}
	\title[]{$p$-Laplacian operator with potential in generalized Morrey Spaces}
\begin{document}
\maketitle
	\begin{abstract}
	We study some basic properties of generalized Morrey spaces $\mathcal{M}^{p,\phi}(\R^{d})$. Also, the problem $-\mbox{div}(|\nabla u|^{p-2}\nabla u)+V|u|^{p-2}u=0$ in $\Omega$, where $\Omega$ is a bounded open set in $\R^d$, and potential $V$ is assumed to be not equivalent to zero and lies in $\mathcal{M}^{p,\phi}(\Omega)$, is studied. Finally, we establish the strong unique continuation for the $p$-Laplace operator in the case $V\in\mathcal{M}^{p,\phi}(\R^d)$.
	\end{abstract}	

\section{Introduction}
The so called Morrey space were introduced in 1938 (see \cite{Morrey27}) in relation to regularity problems of solution to partial differential equations.

\medskip
Although used quite often in problems related to PDE's, the awareness of such spaces among mathematicians is not so widespread as in the case of Lebesgue or Sobolev spaces. Detailed information about Morrey spaces may be found in \cite{Adams, Sawano}.

\medskip
The Morrey space $\mathcal{M}^{p,\lambda}(\Omega)$ is defined as
\[
\mathcal{M}^{p,\lambda}(\Omega)=\left\{f\in L_p(\Omega):\sup_{\substack{\lambda\in\Omega\\r>0}}\frac{1}{r^{\lambda}}\int\limits_{B(x,r)}|f(y)|^{p}\,dy<\infty\right\}
\]
where $1\leq p<\infty$, $\lambda\geq 0$, $\Omega\subset\R^{d}$ is a bounded open set and $B(x,r)$ denotes the open ball centered at $x$ with radius $r$, i.e. $B(x,r)=\{y\in\Omega:|y-x|<r\}$.

\medskip
This is a Banach space with respect to the norm
\[
\|f\|_{\mathcal{M}^{p,\lambda}(\Omega)}=\sup_{\substack{x\in\Omega\\r\in(0,\diam(\Omega))}}\left(\frac{1}{r^{\lambda}}\int\limits_{B(x,r)}|f(y)|^{p}\,dy\right)^{\frac{1}{p}}.
\]
We will explore some properties of spaces $\mathcal{M}^{p,\lambda}(\Omega)$, aiming to establish some relations with the Fefferman's inequality, the Poisson equation, the $p$-Laplacian and the unique continuation principle. The last is fundamental and of independent interest in the theory of partial differential equations. Application regarding the vanishing of solutions are often associated for instance with solvability, stability, geometrical properties of solutions, and so on.

\medskip
We are concerned with the following problem 
\begin{equation}\label{eq1.2}
-\mbox{div}(|\nabla u|^{p-2}\nabla u)+V|u|^{p-2}u=0\quad\mbox{in $\Omega$}
\end{equation}
and the potential $V$ is assumed to be not equivalent to zero and lies in $\mathcal{M}^{p,\phi}(\Omega)$ for $1<p<\infty$ and $\phi$ is an almost decreasing function which satisfies the doubling condition.

\medskip
Specifically, we are interested in studing a family of functions which enjoys the strong unique continuation property, that is, functions besides the possible zero functions which have zero infinite order.

\medskip
By $L_{p,loc}(\Omega)$ we denoted the set of functions $u$ such that 
\[
\int\limits_{K}|u(x)|^{p}\,dx<\infty
\]
for all compact subset $K\subset\Omega$.
\begin{dfn}\label{def1.1}
We say tha a function $u\in L_{p,loc}(\Omega)$ vanishes of infinite order at point $x_0$ if for any natural number $N$ there exists a constant $C_N$ such that
\begin{equation}\label{eq1.3}
\int\limits_{B(x_0,r)}|u(x)|^p\,dx\leq C_Nr^N
\end{equation}
for all $n\in\N$ and for small positive number $r$.
\end{dfn}
\begin{dfn}
	We say that \eqref{eq1.2} has strong unique continuation property if and only if any solution $u$ of \eqref{eq1.2} in $\Omega$ is identically zero in $\Omega$ provided that $u$ vanishes of infinite order that a point $x_0$ in $\Omega$. 
\end{dfn}
There is an extensive literature on unique continuation property. We refer the reader to the work of Zamboni on unique continuation for nonnegative solutions of quasilinear elliptic equation \cite{Zamboni3}, also the work of Jerison-Kenig on unique continuation for Schr\"odinger operator \cite{Jerison1}. The same work is done by Chiarenza and Frasca but for linear elliptic operator in the case where $V\in L^{\frac{n}{2}}$ when $n>2$ \cite{Chiarenza2}.

\medskip

This paper is organized as follows. In section \ref{sec:2}, we study a generalization of Morrey space. Some basic results, such as embeddings, completeness, etc. are established. In section \ref{sec:3}, we shall investigate the relationship between Fefferman's inequality and the Poisson equation. Finally, section \ref{sec:4} is devoted to investigate the strong unique continuation property of the so called $p$-Laplacian. The $p$-Laplacian has been much studied during the last sixty eight years, although its theory is by now rather developed, some challenging problems remain unsolved.

\section{Generalized Morrey Space}\label{sec:2}

Generalized Morrey space $\mathcal{M}^{p,\phi}(\R^{d})$ are equipped with the parameter, that is, $1\leq p<\infty$, and the function $\phi:(0,\infty)\longrightarrow(0,\infty)$. We assume that $\phi$ is in $\G_p$, which is the set of all functions $\phi:(0,\infty)\longrightarrow(0,\infty)$ such that $\phi$ is almost decreasing, that is $r\leq s$ implies $\phi(r)\geq C\phi(s)$ and $t\mapsto t^{\frac{n}{p}}$ is almost increasing, that is $r\leq s$ implies $r^{\frac{n}{p}}\phi(r)\leq Cs^{\frac{n}{p}}\phi(s)$ for some $C>0$. Note that, $\phi\in \G_p$ implies that $\phi$ satisfies the doubling condition, that is, there exists $C>0$ such that 
\[
\frac{1}{C}\leq\frac{\phi(r)}{\phi(s)}\leq C
\]
for every $r$ and $s$ with $\frac{1}{2}\leq \frac{r}{s}\leq 2$.
\begin{dfn}[Generalized Morrey space]\label{def2.1}
The generalized Morrey Space $\mathcal{M}^{p,\phi}(\R^{d})=\mathcal{M}^{p,\phi}$ is defined as
\[
\mathcal{M}^{p,\phi}(\R^{d})=\{f\in L_p(\R^d):\|f\|_{\mathcal{M}^{p,\phi}}<\infty\}
\]
with
\begin{equation}\label{eq2.1}
\|f\|_{\mathcal{M}^{p,\phi}(\R^{d})}=\sup_{B(a,r)}\frac{1}{\phi(r)}\left(\frac{1}{r^n}\int\limits_{B(a,r)}|f(x)|^p\right)
\end{equation}
where $B(a,r)=\{x\in\R^d:|x-a|<r\}$.
\end{dfn}

Observe that if $1\leq p<\infty$ and $\phi(r)=r^{\frac{\lambda-n}{p}}$ for $\lambda>0$, then $\mathcal{M}^{p,\phi}=\mathcal{M}^{p,\lambda}$ which means that we recovered the space defined in \eqref{eq2.1}. Now it is propitious to point out that $\R^{d}$ is equipped with a Borel measure $\mu$ satisfying the growth condition of order $n$, with $0<n\leq d$, that is, there exist a constant $C>0$ such that
\[
\mu(B(a,r))\leq Cr^{n}
\]
for every ball $B(a,r)$ centered at $a\in\R^{d}$ with radius $r>0$.

\medskip
The following lemma indicates that the characteristic function on some balls is contained in generalized Morrey spaces.
\begin{lema}\label{lem2.1}
Let $1\leq p<\infty$ and $\phi\in\G_p$. If $B_0=B(0,r_0)$ for some $r_0>0$. Then, there exist $C>1$ and $B>0$ such that
\[
\frac{B}{\phi(r_0)}\leq\|\chi_{B_0}\|_{\mathcal{M}^{p,\phi}}\leq\frac{C}{\phi(r_0)}.
\]
\end{lema}
\begin{proof}
By definition of $\|\cdot\|_{\mathcal{M}^{p,\phi}}$, we have 
\begin{align*}
	\|\chi_{B_0}\|_{\mathcal{M}^{p,\phi}}=&\sup_{B(a,r)}\frac{1}{\phi(r)}\left(\frac{1}{r^n}\int\limits_{B(a,r)}|\chi_{B_0}(x)|^{p}\,d\mu\right)^{\frac{1}{p}}\\
		\geq&\frac{1}{\phi(r_0)}\left(\frac{\mu(B_0)}{r_{0}^n}\right)^{\frac{1}{p}}\\
		=&\frac{B}{\phi(r_0)}
\end{align*}
where $B=\left(\frac{\mu(B_0)}{r_0^n}\right)^{\frac{1}{p}}$, proving the first inequality.

\medskip
For the second inequality, we separate the proof into two cases:

\medskip
First case $r\leq r_0$, for this case, we have 
\[
\phi(r)\geq C\phi(r_0).
\]
Thus 
\begin{align*}
	\frac{1}{\phi(r)}\left(\frac{1}{r^n}\int\limits_{B(a,r)}|\chi_{B_0}(x)|^{p}\,d\mu\right)^{\frac{1}{p}}\leq&\frac{C}{\phi(r_0)}\left(\frac{\mu(B(a,r)\cap B_0)}{r^n}\right)^{\frac{1}{p}}\\
		\leq&\frac{C}{\phi(r_0)}.
\end{align*}
Second case $r\geq r_0$ Since 
\[
r_0^{\frac{n}{p}}\phi(r_0)\leq Cr^{\frac{n}{p}}\phi(r),
\]
we have
\begin{align*}
	\frac{1}{\phi(r)}\left(\frac{1}{r^n}\int\limits_{B(a,r)}|\chi_{B_0}(x)|^p\,d\mu\right)^{\frac{1}{p}}\leq& Cr^{\frac{n}{p}}r_0^{-\frac{n}{p}}\left(\frac{\mu(B(a,r)\cap B_0)}{r^n}\right)^{\frac{1}{p}}\\
		\leq&\frac{C r^{\frac{n}{p}}r_0^{-\frac{n}{p}}}{\phi(r_0)}\left(\frac{\mu(B_0)}{r^n}\right)^{\frac{1}{p}}\\
		\leq&\frac{C}{\phi(r_0)}.
\end{align*}
From these two cases, we can conclude that
\[
\|\chi_{B_0}\|_{\mathcal{M}^{p,\phi}}\leq \frac{C}{\phi(r_0)}.
\]
Thus we have proved the inequality.
\end{proof}

The following lemma shows us some conditions on the functions $\phi_1,\phi_2$ for which we can compare the spaces $\mathcal{M}^{p_1,\phi_1}$ and $\mathcal{M}^{p_2,\phi_2}$. The result is as follows. 
\begin{lema}
Let $1\leq p_1\leq p_2<\infty$, $\phi_1,\phi_2\in\G_p$. Then the following statements are equivalent:
\begin{enumerate}[\normalfont(I)]
	\item $\phi_2\leq C\phi_1$
	\item $\mathcal{M}^{p_2,\phi_2}\subset \mathcal{M}^{p_1,\phi_1}$ with 
	\[
	\|f\|_{\mathcal{M}^{p_1,\phi_1}}\leq C\|f\|_{\mathcal{M}^{p_2,\phi_2}}
	\]
	for every $f\in\mathcal{M}^{p_2,\phi_2}$.
\end{enumerate}
\end{lema}
\begin{proof}
Assume that (I) holds and let $f\in\mathcal{M}^{p_2,\phi_2}$ hence by H\"older's inequality we have 
\begin{align*}
	&\frac{1}{\phi_1(r)}\left(\frac{1}{r^n}\int\limits_{B(a,r)}|f(x)|^{p_1}\right)^{\frac{1}{p_1}}\\
		&\phantom{\frac{1}{\phi_1(r)}}\leq  \frac{C}{\phi_2(r)}\left[\left(\frac{1}{r^n}\int\limits_{B(a,r)}(|f(x)|^{p_1})^{\frac{p_2}{p_1}}\,d\mu\right)^{\frac{p_1}{p_2}}\left(\int\limits_{B(a,r)}d\mu\right)^{1-\frac{p_1}{p_2}}\right]^{\frac{1}{p_1}}\\
		&\phantom{\frac{1}{\phi_1(r)}}\leq \frac{C}{\phi_2(r)}\left[\left(\frac{1}{r^n}\int\limits_{B(a,r)}|f(x)|^{p_2}\,d\mu\right)^{\frac{p_1}{p_2}}\left(\mu(B(a,r))\right)^{1-\frac{p_1}{p_2}}\right]^{\frac{1}{p_1}}\\
		&\phantom{\frac{1}{\phi_1(r)}}\leq \frac{C}{\phi_2(r)}\left[\left(\frac{1}{r^n}\int\limits_{B(a,r)}|f(x)|^{p_2}\,d\mu\right)^{\frac{p_1}{p_2}}\left(r^n\right)^{1-\frac{p_1}{p_2}}\right]^{\frac{1}{p_1}}\\
		&\phantom{\frac{1}{\phi_1(r)}}\leq \frac{C}{\phi_2(r)}\left(\frac{1}{r^n}\int\limits_{B(a,r)}|f(x)|^{p_2}\,d\mu\right)^{\frac{1}{p_2}}\\
		&\phantom{\frac{1}{\phi_1(r)}}\leq C\|f\|_{\mathcal{M}^{p_2,\phi_2}}.
\end{align*}
By taking the supremum of the left side for every $a$ and $r$ we obtain
\[
\|f\|_{\mathcal{M}^{p_1,\phi_1}}\leq C\|f\|_{\mathcal{M}^{p_2,\phi_2}}.
\]
Now, assume that (II) hods. Let $B_0=B(0,r_0)$ as before, then
\begin{equation}\label{eq2.2}
\|\chi_{B_0}\|_{\mathcal{M}^{p_1,\phi_1}}\leq C\|\chi_{B_0}\|_{\mathcal{M}^{p_2,\phi_2}}.
\end{equation}
By Lemma \ref{lem2.1} we have 
\begin{equation}\label{eq2.3}
\frac{1}{\phi_1(r_0)}\leq \|\chi_{B_0}\|_{\mathcal{M}^{p_1,\phi_1}}
\end{equation}
and
\begin{equation}\label{eq2.4}
	\|\chi_{B_0}\|_{\mathcal{M}^{p_1,\phi_1}}\leq\frac{C}{\phi_2(r_0)}.
\end{equation}
By \eqref{eq2.2}, \eqref{eq2.3} and \eqref{eq2.4} we arrive at
\[
\frac{1}{\phi_1(r_0)}\leq\frac{C}{\phi_2(r_0)},
\]
and so
\[
\phi_2(r_0)\leq C\phi_1(r_0).
\]
Since $r_0$ is any positive real number, we have
\[
\phi_2\leq \phi_1.
\]
\end{proof}

The next result gives us a H\"older's type inequality.
\begin{thm}\label{thm2.1}
Let $f\in\mathcal{M}^{p,\phi}(\R^d)$ and $g\in\mathcal{M}^{q,\phi}(\R^d)$. Then 
\[
\|fg\|_{\mathcal{M}^{1,\phi}}\leq\phi(r)\|f\|_{\mathcal{M}^{p,\phi}}\|g\|_{\mathcal{M}^{q,\phi}}
\]
with $\frac{1}{p}+\frac{1}{q}=1$.
\end{thm}
\begin{proof}
By H\"older's  inequality we have 
\begin{align*}
	&\frac{1}{\phi(r)r^n}\int\limits_{B}|f(x)g(x)|\,d\mu=\frac{\phi(r)}{[\phi(r)]^2}\int\limits_{B}\left|\frac{f(x)}{r^{\frac{n}{p}}}\right|\left|\frac{g(x)}{r^{\frac{n}{q}}}\right|\,d\mu\\
		&\hspace{.5cm}\leq\phi(r)\left(\frac{1}{\phi(r)}\left(\frac{1}{r^n}\int\limits_{B}|f(x)|^p\,d\mu\right)^{\frac{1}{p}}\left(\frac{1}{\phi(r)}\left(\frac{1}{r^n}\int\limits_{B}|g(x)|{ q}\,d\mu\right)^{\frac{1}{q}}\right)\right)\\
		&\hspace{.5cm}\leq\phi(r)\sup_{B}\frac{1}{\phi(r)}\left(\frac{1}{r^n}\int\limits_{B}|f(x)|^{p}\,d\mu\right)^{\frac{1}{p}}\sup_{B}\frac{1}{\phi(r)}\left(\frac{1}{r^n}\int\limits_{B}|g(x)|^{q}\,d\mu\right)^{\frac{1}{q}}\\
		&\hspace{.5cm}=\phi(r)\|f\|_{\mathcal{M}^{p,\phi}}\|g\|_{\mathcal{M}^{q,\phi}}.
\end{align*}
And so 
\[
\sup_{B}\frac{1}{\phi(r)}\left(\frac{1}{r^n}\int\limits_{B}|f(x)g(x)|\,d\mu\right)\leq\phi(r)\|f\|_{\mathcal{M}^{p,\phi}}\|g\|_{\mathcal{M}^{q,\phi}}.
\]
That is 
\[
\|fg\|_{\mathcal{M}^{1,\phi}}\leq \phi(r)\|f\|_{\mathcal{M}^{p,\phi}}\|g\|_{\mathcal{M}^{q,\phi}},
\]
and so the theorem is proved.
\end{proof}
In the coming result we derive the Minkowski inequality, which  in this case is independent of Theorem \ref{thm2.1}.
\begin{thm}\label{thm2.2}
If $f$ and $g$ belong to $\mathcal{M}^{p,\phi}(\R^d)$, then $f+g$ belongs to $\mathcal{M}^{p,\phi}(\R^d)$. Moreover,
\[
\|f+g\|_{\mathcal{M}^{p,\phi}}\leq\|f\|_{\mathcal{M}^{p,\phi}}+\|g\|_{\mathcal{M}^{p,\phi}}.
\]
\end{thm} 
\begin{proof}
Since $f$ and $g$ belong to $\mathcal{M}^{p,\phi}(\R^d)$ by definition $f$ and $g$ belong to $L_p(\R^d)$. Hence the Minkowski inequality holds. So
\[
\left(\int\limits_{B}|f+g|^p\,d\mu\right)^{\frac{1}{p}}\leq\left(\int\limits_{B}|f|^p\,d\mu\right)^{\frac{1}{p}}+\left(\int\limits_{B}|g|^p\,d\mu\right)^{\frac{1}{p}}.
\]	
Thus 
\begin{align*}
	&\frac{1}{\phi(r)r^{\frac{n}{p}}}\left(\int\limits_{B}|f+g|^p\,d\mu\right)^{\frac{1}{p}}\\
	&\phantom{\frac{1}{\phi(r)r^{\frac{n}{p}}}}\leq\frac{1}{\phi(r)r^{\frac{n}{p}}}\left(\int\limits_{B}|f|^p\,d\mu\right)^{\frac{1}{p}}+\frac{1}{\phi(r)r^{\frac{n}{p}}}\left(\int\limits_{B}|g|^p\,d\mu\right)^{\frac{1}{p}}\\
		&\phantom{\frac{1}{\phi(r)r^{\frac{n}{p}}}}=\frac{1}{\phi(r)}\left(\frac{1}{r^n}\int\limits_{B}|f|^p\,d\mu\right)^{\frac{1}{p}}+\frac{1}{\phi(r)}\left(\frac{1}{r^n}\int\limits_{B}|g|^p\,d\mu\right)^{\frac{1}{p}}\\
		&\phantom{\frac{1}{\phi(r)r^{\frac{n}{p}}}}\leq\sup_{B}\frac{1}{\phi(r)}\left(\frac{1}{r^n}\int\limits_{B}|f|^p\,d\mu\right)^{\frac{1}{p}}+\sup_{B}\frac{1}{\phi(r)}\left(\frac{1}{r^n}\int\limits_{B}|g|^p\,d\mu\right)^{\frac{1}{p}}\\
		&\phantom{\frac{1}{\phi(r)r^{\frac{n}{p}}}}=\|f\|_{\mathcal{M}^{p,\phi}}+\|g\|_{\mathcal{M}^{p,\phi}}.
\end{align*}
Finally
\[
\frac{1}{\phi(r)}\left(\frac{1}{r^n}\int\limits_{B}|f+g|^p\,d\mu\right)^{\frac{1}{p}}\leq \|f\|_{\mathcal{M}^{p,\phi}}+\|g\|_{\mathcal{M}^{p,\phi}}
\]
holds for any ball $B$. Therefore 
\[
\sup_{B}\frac{1}{\phi(r)}\left(\frac{1}{r^n}\int\limits_{B}|f+g|^p\,d\mu\right)^{\frac{1}{p}}\leq \|f\|_{\mathcal{M}^{p,\phi}}+\|g\|_{\mathcal{M}^{p,\phi}}
\]
and so
\[
\|f+g\|_{\mathcal{M}^{p,\phi}}\leq \|f\|_{\mathcal{M}^{p,\phi}}+\|g\|_{\mathcal{M}^{p,\phi}}.
\]
\end{proof}

Now, we are ready to prove the completeness of $\mathcal{M}^{p,\phi}$ space. 

\begin{thm}\label{thm2.3}
$\mathcal{M}^{p,\phi}$, equipped with the norm \eqref{eq2.1}, is a Banach space. 
\end{thm}
\begin{proof}
Let $\{f_n\}_{n\in\N}$ be a Cauchy sequence in $\mathcal{M}^{p,\phi}(\R^d)$, since $\mathcal{M}^{p,\phi}(\R^d)\subset L_p(\R^d)$, then $\{f_n\}_{n\in\N}$ is a Cauchy sequence in $L_p(\R^d)$, therefore there exists an $f\in L_p(\R^d)$ such that $f_n\to f$ in $L_p(\R^d)$.

\medskip
Thus for $\epsilon>0$ and $r>0$ there exist $n_r>0$ such that 
\[
\|f-f_n\|_{L_p}<\phi(r)r^{\frac{n}{p}}\epsilon\quad\mbox{if}\quad n\geq n_r,
\]
from this latter inequality we arrived at 
\begin{equation}\label{eq2.5}
\|f-f_n\|_{\mathcal{M}^{p,\phi}}<\epsilon\quad\mbox{if}\quad n\geq n_r.
\end{equation}
Applying Theorem \ref{thm2.2} and \eqref{eq2.5} we get
\begin{align*}
	\|f\|_{\mathcal{M}^{p,\phi}}\leq&\|f-f_n\|_{\mathcal{M}^{p,\phi}}+\|f_{n_r}\|_{\mathcal{M}^{p,\phi}}\\
		<&\epsilon+\|f_{r_n}\|\leq C.
\end{align*}
And so $f\in\mathcal{M}^{p,\phi}(\R^d)$ which end the proof of Theorem \ref{thm2.3}
\end{proof}

We end this section by proving a Hedberg-type inequality. Before doing so, let us recall two well-known definitions.

\medskip
The first one is the maximal operator. Given a function $f\in L_{1,loc}(\R^d)$ the Hardy-Littlewood maximal function for $x\in\R^d$ it is defined as
\begin{equation}\label{eq2.6}
Mf(x)=\sup_{r>0}\frac{1}{\mu(B(x,r))}\int\limits_{B(x,r)}|f(y)|\,d\mu.
\end{equation}
The second one is the so-called Riesz potential operator which is defined as 
\begin{equation}\label{eq2.7}
I_{\alpha}f(x)=\int\limits_{\R^d}\frac{f(y)}{|x-y|^{n-\alpha}}
\end{equation}
for $0\leq\alpha<n\leq d$. Now, we are ready to establish our announced result.
\begin{thm}\label{thm2.4}
Suppose that for some $0\leq \lambda<n-1$ we have 
\[
\int_{r}^{\infty}\phi(t)\,dt\leq Cr^{\lambda+1-n},\quad r>0.
\]
Then, for any $f\in\mathcal{M}^{p,\phi}$, we have the following pointwise inequality
\[
|I_1f(x)|\leq C|Mf(x)|^{1-\frac{1}{n-\lambda}}\|f\|_{\mathcal{M}^{p,\phi}}^{\frac{1}{n-\lambda}}
\]
for $x\in\R^{d}$.
\end{thm}
\begin{proof}
Let $f\in\mathcal{M}^{p,\phi}$ and $x\in\R^d$. For every $r>0$ we have 
\begin{align*}
	|I_1f(x)|\leq&\int\limits_{|x-y|<r}\frac{|f(y)|d\mu}{|x-y|^{n-1}}+\int\limits_{|x-y|\geq r}\frac{|f(y)|d\mu}{|x-y|^{n-1}}\\
		=&A+B.
\end{align*}
Observe that for the first integral we obtain
\begin{align}
	A=&\int\limits_{|x-y|<r}\frac{|f(y)|d\mu}{|x-y|^{n-1}} \notag\\
		=&\sum_{j=-\infty}^{-1}\int\limits_{2^{j}r<|x-y|\leq2^{j+1}r}\frac{|f(y)|d\mu}{|x-y|^{n-1}} \notag\\
		\leq&\sum_{j=-\infty}^{-1}\frac{1}{(2^jr)^{n-1}}\int\limits_{|x-y|\leq2^{j+1}r}|f(y)|\,d\mu \notag\\
		=&\sum_{j=-\infty}^{-1}2^n(2^jr)\frac{1}{(2^{j+1}r)^{n}}\int\limits_{B(x,2^{j+1}r)}|f(y)|\,d\mu \notag\\
		\leq&2^nrMf(x)\sum_{j=-\infty}^{-1}2^j\notag\\
		\leq&CrMf(f).\label{eq2.8}
\end{align}
Meanwhile, for the second integral, we have the following estimate
\begin{align*}
	B=&\int\limits_{|x-y|\geq r}\frac{|f(y)|}{|x-y|^{n-1}}\,d\mu\\
		=&\sum_{j=0}^{\infty}\int\limits_{2^jr<|x-y|\leq2^{j+1}r}\frac{|f(y)|d\mu}{|x-y|^{n-1}}\\
		\leq&\sum_{j=0}^{\infty}\frac{1}{(2^jr)^{n-1}}\int\limits_{|x-y|\leq 2^{j+1}r}|f(y)|\,d\mu\\
		=&\sum_{j=0}^{\infty}2^n(2^{j}r)\frac{1}{(2^{j+1}r)^{n}}\int\limits_{B(x,2^{j+1}r)}|f(y)|\,d\mu\\
		\leq&\sum_{j=0}^{\infty}2^n(2^{j}r)\frac{1}{(2^{j+1}r)^{n}}\left(\int\limits_{B(x,2^{j+1}r)}|f(y)|^{p}\,d\mu\right)^{\frac{1}{p}}\left(\frac{1}{(2^{j+1}r)^{n}}\mu(B(x,2^{j+1}r))\right)^{\frac{1}{q}}\\
		\leq&C\sum_{j=0}^{\infty}2^n(2^{j}r)\|f\|_{\mathcal{M}^{p,\phi}}\left(\frac{1}{(2^{j+1}r)^n}(2^{j+1}r)^n\right)^{\frac{1}{q}}\phi(2^{j+1}r)\\
		=&C\left[\sum_{j=0}^{\infty}2^n(2^{j}r)\phi(2^{j+1}r)\right]\|f\|_{\mathcal{M}^{p,\phi}}.
\end{align*}

Since $\phi$ is almost decreasing, we observe that for $j=0,1,2,\dots$
\[
(2^jr)\phi(2^{j+1}r)\leq C\int_{2^jr}^{2^{j+1}r}\phi(t)\,dt.
\]
this last inequality and our assumption then lead us to 
\begin{align*}
	B\leq& C\|f\|_{\mathcal{M}^{p,\phi}}\sum_{j=0}^{\infty}\int_{2^jr}^{2^{j+1}r}\phi(t)\,dt\\
		\leq&C\|f\|_{\mathcal{M}^{p,\phi}}\int_{r}^{\infty}\phi(t)\,dt\\
		\leq&Cr^{\lambda+1-n}\|f\|_{\mathcal{M}^{p,\phi}}.
\end{align*}
Now, by choosing
\[
r=\left(\frac{Mf(x)}{\|f\|_{\mathcal{M}^{p,\phi}}}\right)^{\frac{1}{\lambda-n}},
\]
we obtain 
\begin{align*}
	|I_1f(x)|\leq& Cr\left(Mf(x)+r^{\lambda-n}\|f\|_{\mathcal{M}^{p,\phi}}\right)\\
		\leq&C[Mf(x)]^{1-\frac{1}{n-\lambda}}\|f\|_{\mathcal{M}^{p,\phi}}^{\frac{1}{\lambda-n}},
\end{align*}
this completes the proof.
\end{proof}

\section{Fefferman inequality and its realtion with the Poisson equation.}\label{sec:3}

Let us consider the following problem (Dirichlet problem)
\begin{equation}\label{eq3.1}
\left\lbrace
\begin{aligned}
	-\Delta z=&V\quad\mbox{on $B$}\\
	z=&0\quad\mbox{on $\partial B$.}
\end{aligned}
\right.
\end{equation}
The equation $-\Delta z=V$ is known as the Poisson equation. It is well known that the solution of the above problem is given by the convolution 
\[
z(x)=\int\limits_{\R^d}\phi(x-y)V(y)\,dy,
\]
where $\phi$ is the fundamental solution of the Laplace equation, and if $V\in C_{c}^{2}(\R^d)$ it is clear that $z\in C_{c}^{2}(\R^d)$, (see \cite{evans2022partial} for details).

\medskip
Furthermore, $z$ may be written as
\[
z(x)=\frac{1}{\omega_{n-1}}\int\limits_{\R^d}\frac{\nabla z(y)\cdot(x-y)}{|x-y|^n}\,dy
\]
where $\omega_{n-1}$ represents the $(n-1)$-dimensional measure of the sphere $S^{n-1}$. Then
\[
\nabla z(x)=\frac{1}{\omega_{n-1}}\int\limits_{\R^d}\frac{\nabla^2 z(y)\cdot(x-y)}{|x-y|^n}\,dy.
\]
Thus 
\[
|\nabla z(x)|=\frac{1}{\omega_{n-1}}\int\limits_{\R^d}\frac{|\nabla^2 z(y)|}{|x-y|^{n-1}}\,dy
\]
it is known that $\nabla^2z=\Delta z$.
\begin{dfn}\label{def3.1}
A function $\omega(x)\geq 1$ is said to be $A_1$-class if 
\[
M\omega(x)\leq C_1\omega(x)
\]
for almost $x\in\R^d$ and for some constant $C_1>0$.
\end{dfn}
Our next task is to state and prove the Fefferman's inequality, assuming that $V$ belongs to
\[
A_1\cap\mathcal{M}^{p,\phi}(\R^d)\cap C_{c}^{2}(\R^d)
\]
with $1\leq p<q\leq\frac{n-\lambda}{n-\lambda-1}$.

\medskip
Fefferman's inequality has been shown in different settings. For instance, in \cite{Zamboni3}, Zamboni proved the Fefferman inequality allowing $V$ to be in a generalized Kato class. See also \cite{Castillo5}.

\medskip
On the other hand, Castillo, Ramos and Rojas proved the Fefferman's inequality allowing $V$ to belong to the Kato class with $p=2$, as well as Castillo, Rafeiro and Rojas proved the Fefferman's inequality allowing $V$ to be in $L_{\frac{n}{p}}(\Omega)$ (see \cite{Castillo-Rafeiro-Rojas*}). For definitions and details, see \cite{Castillo4,Castillo5,Castillo6,Castillo7,Zamboni3}

\begin{thm}[Fefferman's inequality on $\mathcal{M}^{p,\phi}$]\label{thm3.1}
Suppose that for some $0\leq\lambda< n-1$ we have 
\[
\int_{r}^{\infty}\phi(t)\,dt\leq Cr^{\lambda+1-n}\qquad r>0.
\]
Let
\[
1<p\leq q\leq\frac{n-\lambda}{n-\lambda-1}\quad\mbox{and}\quad V\in A_1\cap\mathcal{M}^{p,\phi}\cap C_{c}^{2}(\R^d).
\]
Then
\[
\int\limits_{\R^d}|u(x)|^pV(x)\,d\mu\leq C(n,p)\|V\|_{\mathcal{M}^{p,\phi}}^{\frac{p}{n-\lambda}}\int\limits_{\R^{d}}|\nabla u(x)|^p\,d\mu.
\]
\end{thm}
\begin{proof}
For any $u\in C_{c}^{2}(\R^d)$, let us consider a ball such that $u\in C_{c}^{2}(B)$ and consider the solution $z$ of the Poisson equation 
\[
\left\lbrace
\begin{aligned}
-\Delta z=&V\quad\mbox{on $B$}\\
z=&0\quad\mbox{on $\partial B$.}
\end{aligned}
\right.
\]
Then, by using Theorem \ref{thm2.4} and H\"older's inequality, we have 
\begin{align*}
	&\int\limits_{\R^d}|u(x)|^{p}V(x)\,d\mu=-\int\limits_{B}|u(x)|^p\Delta z(x)\,d\mu(x)\\
		&=\int\limits_{B}\nabla|u(x)|^p\nabla z(z)\,d\mu(x)\\
		&=p\int\limits_{B}|u(x)|^{p-1}|\nabla u(x)|\nabla z(x)\,d\mu(x)\\
		&\leq p\int\limits_{B}|u(x)|^{p-1}|\nabla u(x)||\nabla z(x)|\,d\mu(x)\\
		&\leq\frac{p}{\omega_{n-1}}\int\limits_{B}|u(x)|^{p-1}|\nabla u(x)|\int\limits_{B}\frac{|\nabla^2 z(y)|}{|x-y|^{n-1}}\,d\mu(y)d\mu(x)\\
		&=\frac{p}{\omega_{n-1}}\int\limits_{B}|u(x)|^{p-1}|\nabla u(x)|\int\limits_{B}\frac{|\Delta z(y)|}{|x-y|^{n-1}}\,d\mu(y)d\mu(x)\\
		&=\frac{p}{\omega_{n-1}}\int\limits_{B}|u(x)|^{p-1}|\nabla u(x)|\int\limits_{B}\frac{|V(y)|}{|x-y|^{n-1}}\,d\mu(y)d\mu(x)\\
		&\leq\frac{Cp}{\omega_{n-1}}\int\limits_{B}|u(x)|^{p-1}|\nabla u(x)||MV(x)|^{1-\frac{1}{n-\lambda}}\|V\|_{\mathcal{M}^{p,\phi}
		}^{\frac{1}{n-\lambda}}\,d\mu(x)\\
		&\leq\frac{Cp}{\omega_{n-1}}\|V\|_{\mathcal{M}^{p,\phi}
		}^{\frac{1}{n-\lambda}}\left[\int\limits_{B}|u(x)|^p[MV(x)]^{q\left(1-\frac{1}{n-\lambda}\right)}\,d\mu(x)\right]^{\frac{1}{q}}\left(\int\limits_{B}|\nabla u(x)|^p\,d\mu\right)\\
		&\leq\frac{Cp}{\omega_{n-1}}\|V\|_{\mathcal{M}^{p,\phi}
		}^{\frac{1}{n-\lambda}}\left(\int\limits_{B}|u(x)|^pV(x)\,d\mu(x)\right)^{\frac{1}{q}}\left(\int\limits_{B}|\nabla u(x)|^{p}\right)^{\frac{1}{p}}.
\end{align*}
And so
\[
\int\limits_{B}|u(x)|^pV(x)\,d\mu\leq\left(\frac{Cp}{\omega_{n-1}}\right)^{p}\|V\|_{\mathcal{M}^{p,\phi}}^{\frac{p}{n-\lambda}}\int\limits_{B}|\nabla u(x)|^p\,d\mu.
\]
This completes the proof.
\end{proof}

From Fefferman's inequality we easily derive the following corollary, which can be obtained through a standard partition of the unity.
\begin{coro}\label{coro3.1}
Let $V\in\mathcal{M}^{p,\phi}(\R^d)$ and let $\Omega$ be a bounded subset of $\R^d$, $\supp V\subseteq\Omega$. Then, for any $\sigma>0$, there exists a positive constant $K$ depending on $\sigma$, such that
\[
\int\limits_{\Omega}|u(x)|^pV(x)\,d\mu\leq\sigma\int\limits_{\Omega}|\nabla u(x)|^p\,d\mu+K(\sigma)\int\limits_{\Omega}\int\limits_{\Omega}|u(x)|^p\,d\mu
\]
for all $u\in C_{0}^{\infty}(\Omega)$.
\end{coro}
\begin{proof}
Let $\sigma>0$. Let $r$ be a positive number that will be chosen later. Let $\{\alpha_{k}^{p}\}$, $k=1,2,\dots,N(r)$, be a finite partition of the unity of $\overline{\Omega}$, such that $\supp V\subseteq B(x_k,r)$ with $x_k\in\overline{\Omega}$. We apply Theorem \ref{thm3.1} to the functions $\alpha_k$ and we get
\begin{align}
	&\int\limits_{\Omega}|u(x)|^pV(x)\,d\mu\notag\\
	&=\int\limits_{\Omega}V(x)|u(x)|^p\sum_{k=1}^{N(r)}\alpha_{k}^{p}(x)\,d\mu\notag\\
	&=\sum_{k=1}^{N(r)}\int\limits_{\Omega}V(x)|u(x)\alpha_{k}(x)|^p\,d\mu\notag\\
	&\leq\sum_{k=1}^{N(r)} C\|V\|_{\mathcal{M}^{p,\phi}}^{\frac{p}{n-\lambda}}\left(\int\limits_{\Omega}|\nabla u(x)|^p\alpha_{k}^{p}(x)\,d\mu+\int\limits_{\Omega}|\nabla u(x)|^p||u(x)|^p\,d\mu\right)\notag\\
	&\leq C\|V\|_{\mathcal{M}^{p,\phi}}^{\frac{p}{n-\lambda}}\left(\int\limits_{\Omega}|\nabla u(x)|^p\,d\mu+\frac{N(r)}{r^p}\int\limits_{\Omega}|u(x)|^p\,d\mu\right)\label{eq3.2}
\end{align}
Finally, to obtain the result, it is sufficient to choose $r$ such that
\[
C\|V\|_{\mathcal{M}^{p,\phi}}^{\frac{p}{n-\lambda}}=\sigma.
\]
After that, we note that $N(r)\sim r^{-n}$ and the corollary follows.
\end{proof}
\begin{lema}\label{lem3.1}
Let $B_r$ and $B_{2r}$ be two concentric balls contained in $\Omega$. Then 
\begin{equation}\label{eq3.3}
\int\limits_{B_r}|\nabla u(x)|^{p}\,d\mu\leq\frac{C}{r^p}\int\limits_{B_{2r}}|u(x)|^{p}\,d\mu,
\end{equation}
where the constant $C$ does not depend on $r$, and $u\in C_{0}^{\infty}(B_r)$.
\end{lema}
\begin{proof}
Take $\varphi\in C_{0}^{\infty}(\Omega)$, with $\supp\varphi\subset B_{2r}$, $\varphi(x)=1$ for $x\in B_r$, and $|\nabla\varphi|\leq\frac{C}{r}$. By using $\varphi^{p}$ as a test function in \eqref{eq1.2}; we have
\begin{equation}\label{eq3.4}
\int\limits_{B_{2r}}-\mbox{div}(|\nabla u|^{p-2}\nabla u)\varphi^{p}(u)+\int\limits_{B_{2r}}V|u|^{p-2}u\varphi^pu=0.
\end{equation}
Thus 
\begin{equation}\label{eq3.5}
\int\limits_{B_{2r}}|\nabla u|^{p}\varphi^{p}=-p\int\limits_{B_{2r}}|\nabla u|^{p-2}\varphi^{p-2}\nabla u\cdot\nabla\varphi(\varphi u)-\int\limits_{B_{2r}}V|\varphi u|^{p}.
\end{equation}
With the help of Young's inequality for $\frac{p-1}{p}+\frac{1}{p}=1$, we can estimate the first integral in the right-hand side of \eqref{eq3.5} by 
\begin{equation}\label{eq3.6}
(p-1)\epsilon^{\frac{p}{p-1}}\int\limits_{B_{2r}}|\nabla u|^p\varphi^p+\epsilon^{-p}\int\limits_{B_{2r}}|\nabla\varphi|^{p}|u|^{p}.
\end{equation}
Also, by Corollay \ref{coro3.1} we can estimate the second integral in the right-hand side of \eqref{eq3.5} by
\begin{equation}\label{eq3.7}
\epsilon\int\limits_{B_{2r}}|\nabla(\varphi u)|^p+C_{\epsilon}\int\limits_{B_{2r}}|\varphi u|^p.
\end{equation}
Using this estimates in \eqref{eq3.5} we arrive at
\begin{align*}
	\int\limits_{B_{2r}}|\nabla u|^{p}\varphi^p\leq&\left((o-1)\epsilon^{\frac{p}{(p-1)}}+\epsilon\right)\int\limits_{B_{2r}}|\nabla u|^{p}\varphi^p\\
		&+(\epsilon^p+\epsilon)\int\limits_{B_{2r}}|u|^{p}|\nabla u|^{p}+C_{\epsilon}\int\limits_{B_{2r}}|\nabla u|^p|\varphi|^p.
\end{align*}
Using the fact that $|\nabla\varphi|\leq\frac{C}{v}$, $|\varphi|\leq\frac{C}{v}$ and $\varphi=1$ in $B_{r}$, we immediately get inequality \eqref{eq3.3}.
\end{proof}
\begin{lema}\label{lem3.2}
Let $u\in C_{0}^{\infty}(\Omega)$ where $B_r$ is the ball of radius $r$ in $\R^d$ and let $E=\{x\in B_r:u(x)=0\}$. Then there exits a constant $\beta$ (depending only on $n$) such that
\begin{equation}\label{eq3.8}
\int\limits_{A}|u|\leq\beta\frac{r^n}{\mu(E)}\left(\me(E)\right)^\frac{1}{n}\int\limits_{B_{2r}}|\nabla u|,
\end{equation}
for all ball $B_r$, $u$ as above and all measurable sets $A\subset B_r$.
\end{lema}
\begin{proof}
Let us begin by choosing $r$ such that $\me(A)=\me(B)$ where $m$ is the Lebesgue measure, next, for $x\in B_r$ and $t\in E$ we obtain 
\begin{equation}\label{eq3.9}
-u(x)=u(t)-u(x)=\int_{0}^{|x-y|}\frac{du(x-s\omega)}{d\mu(s)}\,d\mu(s)
\end{equation}
where $\omega=\frac{y-x}{|y-x|}$. Now let us integrate equation \eqref{eq3.9} with respect to $t\in E$, once it is done, we have 
\begin{equation}\label{eq3.10}
-\mu(E)u(x)=\int\limits_{E}du(t)\int_{0}^{|x-y|}\frac{du}{d\mu(s)}(x+s\omega)\,d\mu(s).
\end{equation}
Now, we need to find a bound for \eqref{eq3.10}. By using the Fubini's Theorem,
\begin{align*}
	&\left|-\mu(E)u(x)\right|=\left|\int\limits_{E}d\mu(t)\int_{0}^{|x-y|}\frac{du}{d\mu(s)}(x+s\omega)\,d\mu(s)\right|\\
		&\leq\int\limits_{E}du(t)\int_{0}^{|x-y|}\left|\frac{du}{d\mu(s)}(x+s\omega)\right|\,d\mu(s)\\
		&\leq\int\limits_{B_r}du(t)\int_{0}^{|x-y|}\left|\frac{du}{d\mu(s)}(x+s\omega)\right|\,d\mu(s)\\
		&\leq\int\limits_{B_{2r}}du(t)\int_{0}^{|x-y|}\left|\frac{du}{d\mu(s)}(x+s\omega)\right|\,d\mu(s)\\
		&=\mu(B_{2r})\int_{0}^{|x-y|}\left|\frac{du}{d\mu(s)}(x+s\omega)\right|\,d\mu(s)\\
		&\leq\frac{C(2r)^n}{\me(B(0,1))}\me(B(0,1))\int_{0}^{|x-y|}\left|\frac{du}{d\mu(s)}(x+s\omega)\right|\,d\mu(s)\\
		&=\frac{C(2r)^nn\me(B(0,1))}{n\me(B(0,1))}\int_{0}^{|x-y|}\left|\frac{du}{d\mu(s)}(x+s\omega)\right|\,d\mu(s)\\
		&=\frac{C}{\me(\partial B(0,1))}\int\limits_{\partial B(0,1)}d\omega\int_{0}^{2r}t^{n-1}\,dt\int_{0}^{|x-y|}\left|\frac{du}{d\mu(s)}(x+s\omega)\right|\,d\mu(s)\\
		&=\frac{C}{\me(\partial B(0,1))}\int_{0}^{2r}t^{n-1}\,dt\int\limits_{\partial B(0,1)}d\omega\int_{0}^{|x-y|}\frac{1}{s^{n-1}}\left|\frac{du}{d\mu(s)}(x+s\omega)\right|s^n\,d\mu(s)\\
		&\leq\frac{C}{\me(\partial B(0,1))}\int_{0}^{2r}t^{n-1}\,dt\int\limits_{\partial B(0,1)}d\omega\int_{0}^{r}\frac{|\nabla u(y)|}{|x-y|^{n-1}}s^{n-1}\,d\mu(s)\\
		&=\frac{C(2r)^n}{n\me(\partial B(0,1))}\int\limits_{\partial B(0,1)}\int_{0}^{r}\frac{|\nabla u(y)|}{|x-y|^{n-1}}s^{n-1}\,d\mu(s)d\omega\\
		&=\frac{C(2r)^n}{n\me(\partial B(0,1))}\int\limits_{B(x,r)}\frac{|\nabla u(y)|}{|x-y|^{n-1}}s^{n-1}\,d\mu(y).
\end{align*}
Taking into account that we just write $y=x+r\omega$, thus $r\omega=y-x$ and $r=|x-y|$ since $|\omega|=1$. Up to now, we have obtained that 
\begin{equation}\label{eq3.11}
\mu(E)u(x)\leq\frac{C_1(2r^n)}{n}\int\limits_{B(x,r)}\frac{|\nabla u(y)|}{|x-y|^{n-1}}\,d\mu(y)
\end{equation}
where $C_1=\frac{C}{\me(\partial B(0,1))}$.

\medskip
Hence, let us integrate both sides of \eqref{eq3.11} with respect to $x$ on the set $A$. Indeed
\begin{align*}
	&\mu(E)\int\limits_{A}|u(x)|\,d\mu(x)\leq\frac{C(2r)^n}{n}\int\limits_{A}\int\limits_{B(x,r)}\frac{|\nabla u(y)|}{|x-y|^{n-1}}\,d\mu(y)d\mu(A)\\
		&\phantom{\mu(E)}\leq\frac{C(2r)^n}{n}\int\limits_{B(x,r)}\frac{d\mu(x)}{|x-y|^{n-1}}\int\limits_{B(x,r)}|\nabla u(y)|d\mu(y)\\
		&\phantom{\mu(E)}=\frac{C(2r)^n}{n}n\me(B(0,1))\int_{0}^{r}t^{1-r}t^{n-1}\,dt\int\limits_{B(x,r)}|\nabla u(y)|d\mu(y)\\
		&\phantom{\mu(E)}=C(2r^n)\me(B(0,1))r\int\limits_{B(x,r)}|\nabla u(y)|\,d\mu(y)\\
		&\phantom{\mu(E)}=C(2r)^{n}(\me(B(0,1)))^{1-\frac{1}{n}}(r^n\me(B(0,1)))^{\frac{1}{n}}\int\limits_{B(x,r)}|\nabla u(y)|d\mu(y)\\
		&\phantom{\mu(E)}=C(2r)^{n}(\me(B(0,1)))^{1-\frac{1}{n}}(r^n\me(B))^{\frac{1}{n}}.
\end{align*}
Finally
\[
\int\limits_{A}|u(x)|\,d\mu(x)\leq\frac{\beta r^n(\me(A))^{\frac{1}{n}}}{\mu(E)}\int\limits_{B(x,2r)}|\nabla u(y)|\,d\mu(y)
\]
where
\[
\beta=2^n(\me(B(0,1)))^{1-\frac{1}{n}}.
\]
\end{proof}
\section{Strong Unique Continuation}\label{sec:4}
In this section, we proceed to establish the strong unique continuation property of the eigenfunction for the $p$-Laplacian operator in the case $V\in\mathcal{M}^{p,\phi}(\R^d)$.
\begin{thm}\label{thm4.1}
Let $u\in C_0^{\infty}(\Omega)$ be a solution of $\eqref{eq1.2}$. If $u=0$ on a set $E$ of positive measure, then $u$ has zero of infinite orden in $p$-mean.
\end{thm}
\begin{proof}
Let $E=\{x\in B_r:u(x)=0\}$ with $0<r<1$ and let $E^c$ the complement of $E$.

Taking $r_0$ smaller if necessary, we can assume that $B_{r_0}\subset\Omega$. Since $u=0$ on $E$ by Lemma \ref{eq3.2} we have 
\begin{align*}
	\int\limits_{B_r}|u|^{p}=&\int\limits_{B_r\cap E^{c}}|u|^p\leq \frac{\beta r^n(\me(E^c\cap B_r))}{\mu(E\cap B_r)}\int\limits_{B_r}|\nabla|u|^p|\\
		=&C_n\int\limits_{B_r}|u|^{p-1}|\nabla u|
\end{align*}
where
\[
C_n=\frac{p\beta(\me(B(0,1)))^{\frac{1}{n}}}{\mu(E)}.
\]
By H\"older's inequality
\begin{equation}\label{eq4.11}
\int\limits_{B_r}|u|^p\leq C_n\left(\int\limits_{B_r}|\nabla u|^p\right)^{\frac{1}{p}}\left(\int\limits_{B_r}|u|^p\right)^{\frac{p-1}{p}}
\end{equation}
and by using Young's inequality, we get
\begin{equation}\label{eq4.12}
\int\limits_{B_r}|u|^{p}\leq C_n^r\left(r^{p-1}\int\limits_{B_r}|\nabla u|^p+\frac{p-1}{r}\int\limits_{B_r}|u|^{p}\right).
\end{equation}
Finally by Lemma \eqref{eq3.1}, we have
\begin{equation}\label{eq4.13}
\int\limits_{B_r}|u|^p\leq C_n\int\limits_{B_{2r}}|u|^p.
\end{equation}
Now, let us introduce the following function
\begin{equation}\label{eq4.14}
f(r)=\int\limits_{B_r}|u|^p.
\end{equation}
Observe that consequently $r_0$ depends on $n$. Then \eqref{eq4.13} can be written as 
\begin{equation}\label{eq4.15}
f(r)\leq M_n2^{-n}f(2r)\quad\mbox{for $r\leq r_0$}
\end{equation}
where $M_n=2^{-n}C_n$.

\medskip
Integrating \eqref{eq4.15}, we get
\begin{equation}\label{eq4.16}
-f(\rho)\leq M_n2^{-kn}f(2^n\rho)\quad\mbox{if $2^{k-1}\rho\leq r_0$}
\end{equation}
Now given that $0<r<r_0(n)$, choose $k\in\N$ such that 
\begin{equation}\label{eq4.17}
2^{-k}r_0\leq r\leq 2^{-k+1}r_0.
\end{equation}
From \eqref{eq4.16}, we obtain
\begin{equation}\label{eq4.18}
f(r)\leq M_n2^{-kn}f(2^kr)\leq M_n2^{-kn}f(2r).
\end{equation}
Since $2^{-k}\leq\frac{r}{r_0}$, we finally obtain
\begin{equation}\label{eq4.19}
f(r)\leq M_n\left(\frac{r}{r_0}\right)^{n}f(2r_0).
\end{equation}
And thus, we have 
\begin{equation}\label{eq4.20}
\int\limits_{B_r(x_0)}|u(y)|^p\,d\mu(y)\leq M_n\left(\frac{r}{r_0}\right)^nf(2r_0)
\end{equation}
and this shows that \eqref{eq1.3} holds, which means that $u$ has zero of infinite order in $p$-mean at $x_0$.
\end{proof}
\begin{coro}\label{coro4.1}
	Equation \eqref{eq1.2} has strong unique continuation property.
\end{coro}

\rule{0.25\textwidth}{1pt}
\end{document}